\documentclass[11pt]{article}

\usepackage{latexsym}
\usepackage{amsmath,amsfonts,amssymb}
\usepackage{color,epsfig}

\def\leq{\leqslant}
\def\geq{\geqslant}

\numberwithin{equation}{section}

\newcommand{\D}[1]{{\mathbb#1}}
\newcommand{\NN}{{\D{N}}}
\newcommand{\RR}{{\D{R}}}

\newcommand{\ZZ}{{\D{Z}}}

\begin{document}

\title{Dynamics of a $3$ cluster cell-cycle system with positive linear feedback.}
\author{Bastien Fernandez\footnotemark[1] \ and Todd R.\ Young\footnotemark[2]}
\date{}

\maketitle

\begin{abstract}
In this technical note we calculate the dynamics of a linear feedback model of progression in the
cell cycle in the case that the cells are organized into $k=3$ clusters. We examine the
dynamics in detail for a specific subset of parameters with non-empty interior. 

There is an interior fixed point of the Poincar\'{e} map defined by the system. This fixed point
corresponds to a periodic solution with period $T$ in which the three clusters exchange positions
after time $T/3$. We call this solution $3$-cyclic. In all the parameters studied, the fixed point is either:
\begin{itemize}
\item isolated and locally unstable, or,
\item contained in a neutrally stable set of period $3$ points.
\end{itemize}
In the later case the edges of the neutrally stable set are unstable. This case exists if either the
three clusters are isolated from each other, or, if they interact in a non-essential way.
In both cases the orbits of all other interior points are asymptotic to the boundary. Thus $3$-cyclic 
solutions  are practically unstable in the sense the arbitrarily small perturbations may lead to loss of 
stability and eventual  merger of clusters. Since the single cluster solution (synchronization) is the only solution that
 is asymptotically stable, it would seem to be the most likely to be observed in application if
 the feedback is similar to the form we propose and is positive. 
\end{abstract}

\leftline{\small\today.}

\section{Introduction}

\footnotetext[1]{Centre de Physique Th\'{e}orique (UMR 6207 CNRS - Universit\'e Aix-Marseille~II - Universit\'e Aix-Marseille~I - Universit\'e Sud~Toulon-Var) CNRS Luminy Case 907, 13288 Marseille CEDEX 9, France}
\footnotetext[2] {Corresponding author. 321 Morton Hall, Mathematics, Ohio University, Athens, OH 45701, USA,  {\tt young@math.ohiou.edu}}

This technical note is a supplement to the manuscript \cite{YFBMB}. In that work the dynamics
of cells in a cell cycle were considered when cells are subjected to feedback from other
cells in the cell cycle and the response to the feedback is also dependent on the cells location
within the cell cycle. We refer the reader to that manuscript for the biological motivation for
this study and for background on the model we study.

We consider a piecewise affine model for the dynamics of cell populations along the cell cycle. 
Let a population of $N$ cells  be organized in $k$ equal clusters ($k$ divides $N$) labeled by a discrete index 
$i\in\{0,\dots,k-1\}$. For the cluster $i$, 
the progression along the cycle is represented by the periodic variable $x_i\in\RR/\ZZ$ which evolves 
in time according to the entire population status. More specifically, the variables follows the 
deterministic flow associated with a real vector field $f$, {\sl i.e.} $
\dot{x}_i=f_i(,x_i,\bar{x},r,s),\quad i=0,\dots, k-1 $, that is the same for each $i$.
Here $0<s\leq r<1$ are two parameters respectively governing the length of the {\it signaling} 
$S=[0,s)+\ZZ$ and of the {\it responsive} $R=[r,1)+\ZZ$ regions. 
Let $\sigma$ represent the fraction of cells in the signaling region $S$,  {\sl i.e.}
\[
       \sigma = \frac{  \#\{j\ :\ x_j\in [0,s)\ \text{mod}\ 1\} }{k}.
\label{eqn:I}
\]
Then we will consider the case that  $f$ assumes a simple expression
\[
\frac{dx_i}{dt} = f_i(\bar{x},r,s)=\left\{\begin{array}{ccl}
1                   &    \text{if}   &x_i-\lfloor x_i\rfloor\in [0,s)\\
1+ \sigma   &    \text{if}   &x_i-\lfloor x_i\rfloor \in [r,1) ,
\end{array}\right.
\]
where $\lfloor \cdot\rfloor$ denotes the floor function.
In short terms, clusters in the responsive region are accelerated by a fraction equal to the proportion of cells 
in the signaling region. In \cite{YFBMB} we consider systems in which the feedback experienced
by cells $R$ is a general monotone increasing or decreasing function of $\sigma$.

Relabeling, integer translation of coordinates, and time translation are symmetries of the 
dynamics. Thus, one can assume that all coordinates $x_i(0)$ are initially well-ordered and belong 
to the same unit interval, {\it i.e.}\ we have
\[
0=x_0(0)\leq x_i(0)\leq x_{k-1}(0)<1,\quad i=1,\dots,k-2.
\]
The definition of the vector field $f$ implies that the coordinates cannot cross each other as time evolves; 
thus this ordering is preserved under the dynamics. Moreover, the first coordinate $x_0$ must 
eventually reach 1 (not later than at time 1), {\it i.e.}\ there exists $t_\text{R}\leq 1$ such that 
$x_0(t_\text{R})=1$, and more generally, it must reach any positive integer as time runs. 
Thus the set $x_0\in\NN$ defines a Poincar\'e section for the dynamics and the mapping
\[
(x_1(0),x_2(0),\dots,x_{k-1}(0))  \mapsto  (x_1(t_\text{R}),x_2(t_\text{R}),\dots,x_{k-1}(t_\text{R}))
\]
defines the corresponding return map.

We rely on the following considerations. Starting from $t=0$, compute the 
time $t_1$ that $x_{k-1}$ needs to reach 1 and compute the location of the remaining cells 
at this date. Define $F$ to be this mapping. (We assume first that $x_{k-2}(t_1)<1$ for simplicity, but
we will relax this soon.)
\[
F: (x_1(0),x_2(0),\dots,x_{k-1}(0)) \mapsto  (x_0(t_1),x_1(t_1),\dots,x_{k-2}(t_1))
\]
Notice that $x_0(t_1)=t_1$ by assumption on $x_0(0)$.

\begin{figure}[h]
\medskip
\centerline{
\hbox{\psfig{figure=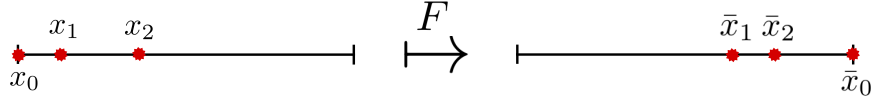,width=11.5cm}}
}
\caption{ Illustration of the map $F$ for $k=3$. Here $F(x_1,x_2) = (x'_1,x'_2) = (x_0(t_1),x_1(t_1))$.}
\label{CASE1}
\end{figure}

Now the time $t_1+t_2$ that $x_{k-2}$ needs to reach 1, together with the population 
configuration at $t=t_1+t_2$, follow by applying $F$ to the configuration 
$(x_0(t_1),x_1(t_1),\dots,x_{k-2}(t_1))$. By repeating the argument,  the desired 
return time $t_\text{R}$ is given by $t_\text{R}=t_1+t_2+\dots +t_k$ and the desired return 
map is $F^k$. Therefore, to understand the dynamics, one only has to compute the first map $F$.

\subsection{General properties for arbitrary $k$}

It was noted in \cite{YFBMB}  that we  may regard $F$ as a continuous piecewise 
affine map of the $(k-1)$-dimensional simplex
$$
0 \leq x_1 \leq x_2 \leq \dots \leq x_{k-1} \leq 1
$$
into itself. (Although the boundaries 0 and 1 are identified in the original flow, in the
analysis here, we consider them as being distinct points for $F$.)

On the edges of the simplex, $F$ has a relatively simple dynamics. Indeed if,
initially, all coordinates are equal, then they must all reach the boundary 1
simultaneously. In other words, on the diagonal ($x_i=x$ for all $i$), we have
$F(x,\dots,x)=(t_1,1,\dots,1)$ where $t_1$ depends on $r,s$ and $x$ (for $x=0$,
we have $t_1=1$ independently of $r$ and $s$). Moreover, starting with $x_{k-1}=1$
implies $t_1=0$ which yields
$$
F(x_1,\dots,x_{k-2},1)=(0,x_1,\dots,x_{k-2})
$$
whatever the remaining coordinates $x_1,\dots,x_{k-2}$ are. As a consequence, the edge
$$
\left\{(x,1,\dots,1)\ :\ x\in [0,1]\right\}
$$
is mapped onto
$$
\left\{(0,x,1,\dots,1)\ :\ x\in [0,1]\right\}
$$
which is mapped onto $\left\{(0,0,x,1,\dots,1)\ :\ x\in [0,1]\right\}$ and so on,
until it reaches the edge $(0,\dots,0,x)$ which is mapped back onto the diagonal (after $k$ iterations).

A particular orbit on the edges is the $k$-periodic orbit passing the vertices, and which
corresponds to the single cluster of velocity 1 in the original flow, namely
\[
(0,\dots,0)\mapsto (1,\dots,1)\mapsto (0,1,\dots,1)\mapsto (0,0,1,\dots,1)\mapsto\dots\mapsto (0,\dots,0,1)\mapsto (0,\dots,0)
\]
From a direct analysis in the original system \cite{YFBMB}, we know that this orbit must be asymptotically
stable for positive feedback and unstable for negative. Since there are two points of a  periodic
orbit at each boundary of every edge (which are themselves globally $k$-periodic 1-dimensional sets),
by the intermediate value theorem, there must be at least one other $k$-periodic orbit on the
edges with coordinates comprised between 0 and 1. Whether this orbit is unique might depend on parameters.

Finally, since the simplex is a convex and compact invariant set under $F$ and the boundary cannot
contain any fixed point, the Brouwer fixed point theorem implies the existence of a fixed point
in its interior.  From the definition of $F$ this corresponds to a solution satrisfying:
\begin{equation}\label{kcyclic}
     x_i(t_1)=x_{i+1}(0) \quad \textrm{for all} \quad i=0, \ldots, k-2,
                 \quad \textrm{and} \quad x_{k-1}(t_1)=x_0(0) \mod 1.
\end{equation}
We call this type of solution a $k$-cyclic solution \cite{YFBMB}.


\section{Dynamics for $k=3$}

We will study the dynamics of the linear model only for a limited subset of parameter space.
It will become clear to the reader that the analysis could be reproduced for the entire
parameter space, but doing so would require a prohibitive amount of time and space.

\subsection{Computation of the map $F$}

When $k=3$, $F$ is defined in the triangle $0\leq x_1\leq x_2\leq 1$ and maps the pair
$(x_1,x_2)$ into $(t_1,x_1(t_1))$, {\it i.e.}\ not only the hitting time $t_1$ needs to
be computed now, but also the location of $x_1$ at this instant has to be specified. There
are numerous cases to study depending on $x_1$ and $x_2$. To list them, we primarily use
the location of $x_2$.
\begin{itemize}
\item[$\bullet$] $0\leq x_2< r-s$. As for $k=2$, in this case $x_2$ moves with velocity 1,
even when lying in $R$. Therefore $t_1=1-x_2$. Since $x_1\leq x_2$, we also have
$x_1\leq r-s$ and thus $x_1(t_1)=x_1+t_1=1-x_2+x_1$.
\item[$\bullet$] $r-s\leq x_2< r$. In this case, $x_2$ enters the responsive region, but only
after a while, precisely at $t=r-x_2$. Its acceleration depends on the location of $x_1$ at this
date, whether $x_1$ belongs to $S$ or not. This alternative is given by the relative sizes of
$r-x_2$ and $s-x_1$. Considering also the possible cases of acceleration for $x_1$ we get
the following sub-cases. {\bf For simplicity, we assume that $s<r-s$ and $r+\frac{5}{3}s<1$}.
\begin{itemize}
\item[$\ast$] $0\leq x_1\leq s$ and $0< r-x_2<s-x_1$. In this case, $x_1$ is still in $S$
when $x_2$ enters $R$. The assumption $r+\frac{5}{3}s<1$ implies that both $x_0$ and $x_1$
get out of $S$ before $x_2$ reaches 1 (existence of the third line in the expression below).
We then have
\[
x_2(t)=\left\{\begin{array}{ccl}
r+\frac{5}{3}(t-r+x_2)&\text{if}&r-x_2< t\leq s-x_1\\
r+\frac{5}{3}(s-x_1-r+x_2)+\frac{4}{3}(t-s+x_1)&\text{if}&s-x_1<t\leq s\\
r+\frac{5}{3}(s-x_1-r+x_2)+\frac{4}{3}x_1+t-s&\text{if}&s<t
\end{array}\right.
\]
which implies $t_1=1-\frac{5}{3}x_2+\frac{x_1}{3}+\frac{2}{3}(r-s)$. Furthermore, the assumption
$s<r-s$ forces $x_0(t)$ to be out of $S$ when $x_1(t)$ enters $R$ (if it does). Then
$x_1(t_1)=x_1+t_1$ in this case.
\item[$\ast$] $0\leq x_1\leq s$ and $s-x_1\leq r-x_2\leq s$ or $s<x_1\leq r-s$ and $x_2< r$.
Here, $x_1(t)$ is out of $S$ when $x_2(t)$ enters $R$. From parameter assumptions we have
$r+\frac{4}{3}s<1$ and thus $x_0$ leaves $S$ before $x_2$ reaches 1. Consequently, the evolution
of $x_2$ reads
\[
x_2(t)=\left\{\begin{array}{ccl}
r+\frac{4}{3}(t-r+x_2)&\text{if}&r-x_2< t\leq s\\
r+\frac{4}{3}(s-r+x_2)+t-s&\text{if}&s\leq t
\end{array}\right.
\]
from which we get $t_1=1-\frac{4}{3}x_2+\frac{r-s}{3}$ and again $x_1(t_1)=x_1+t_1$.
\item[$\ast$] $r-s<x_1\leq x_2$ and $r-s<x_2\leq r$. In this case, $t_1$ remains as before.
However, $x_1$ now suffers some acceleration and we have
\begin{equation}
x_1(t)=\left\{\begin{array}{ccl}
r+\frac{4}{3}(t-r+x_1)&\text{if}&r-x_1< t\leq s\\
r+\frac{4}{3}(s-r+x_1)+t-s&\text{if}&s< t
\end{array}\right.
\label{X1DYN}
\end{equation}
Consequently $x_1(t_1)=1-\frac{4}{3}(x_2-x_1)$.
\end{itemize}
\item[$\bullet$] $r\leq x_2$ and $0\leq x_1\leq s$. For $x_2$ larger than $r$ it is
useful to consider separately the cases $x_1\leq s$ and $x_1>s$. In the latter case,
we always have $x_1(t_1)=x_1+t_1$ but the expression of $t_1$ depends on $x_2$.
We have 3 sub-cases
        \begin{itemize}
        \item[$\ast$] $0\leq x_1\leq s$ and $r\leq x_2\leq 1-\frac{5}{3}s+\frac{x_1}{3}$.
              Here $x_2$ starts sufficiently near 0 so that both $x_1$ and $x_0$ get out of $S$ before $t_1$
        \[
        x_2(t)=\left\{\begin{array}{ccl}
        x_2+\frac{5}{3}t&\text{if}&0< t\leq s-x_1\\
        x_2+\frac{5}{3}(s-x_1)+\frac{4}{3}(t-s+x_1)&\text{if}&s-x_1<t\leq s\\
        x_2+\frac{5}{3}s-\frac{x_1}{3}+t-s&\text{if}&s<t
        \end{array}\right.
        \]
        which implies $t_1=1-x_2+\frac{x_1}{3}-\frac{2}{3}s$.
        \item[$\ast$] $0\leq x_1\leq s$ and $1-\frac{5}{3}s+\frac{x_1}{3}<x_2\leq 1-\frac{5}{3}(s-x_1)$.
              In this case, we have $t_1\leq s$ and hence $x_2$ reaches 1 during the second phase in
              the previous expression. This yields
              $t_1=\frac{3}{4}\left(1-x_2+\frac{x_1}{3}-\frac{s}{3}\right)$.
        \item[$\ast$] $0\leq x_1\leq s$ and $1-\frac{5}{3}(s-x_1)\leq x_2\leq 1$. There $x_2$
              reaches 1 even before $x_1$ reaches s. This results in $t_1=\frac{3}{5}(1-x_2)$.
        \end{itemize}
\item[$\bullet$] $s<x_1$. When $x_1>s$, the situation for $x_2$ is similar to as for $k=2$.
This point may reach 1 before or after $x_0$ leaves $S$. The first case occurs when
$r\leq x_2\leq 1-\frac{4}{3}s$ and gives $t_1=1-x_2-\frac{s}{3}$. In the second case
$1-\frac{4}{3}s<x_2\leq 1$, we have $t_1=\frac{3}{4}(1-x_2)$.
\end{itemize}

We can now focus on $x_1(t_1)$. When $r\leq x_2< 1-\frac{4}{3}s$ and $x_1\leq r-s$, there cannot
be any acceleration and thus $x_1(t_1)=x_1+t_1$. Furthermore, when $1-\frac{4}{3}s\leq x_2\leq 1$,
we have $t_1<s$ and thus $x_1$ receives any acceleration only if it is initially larger than $r-t_1$.
In practice, this leads the following cases.
\begin{itemize}
\item[$\bullet$] $r-s<x_1< r$ and $r\leq x_2< 1-\frac{4}{3}s$. In this case, $s<t_1$ and the
evolution of $x_1$ is as in equation (\ref{X1DYN}). It results that
$x_1(t_1)=1-x_2+\frac{4}{3}x_1-\frac{r}{3}$.
\item[$\bullet$] $r-\frac{3}{4}(1-x_2)<x_1< r$ and $1-\frac{4}{3}s\leq x_2\leq 1$, then only
the first phase in (\ref{X1DYN}) applies and we get again $x_1(t_1)=1-x_2+\frac{4}{3}x_1-\frac{r}{3}$.
\item[$\bullet$] $r\leq x_1<x_2$ and $r\leq x_2< 1-\frac{4}{3}s$. Now, we have
\[
x_1(t)=\left\{\begin{array}{ccl}
x_1+\frac{4}{3}t&\text{if}&t\leq s\\
x_1+\frac{4}{3}s+t-s&\text{if}&s\leq t
\end{array}\right.
\]
Thus $x_1(t_1)=x_1+t_1+\frac{s}{3}=1-x_2+x_1$.
\item[$\bullet$] $r\leq x_1<x_2$ and $1-\frac{4}{3}s<x_2\leq 1$. As $t_1<s$, only the first phase applies in the previous expression. We get $x_1(t_1)=1-x_2+x_1$.
\end{itemize}

All together, when the parameters are such that $2s<r<1-\frac{5}{3}s$, the triangle decomposes
into 13 polygonal subdomains with affine dynamics, see Figure \ref{CASE1}. For simplicity, we have labeled
these domains by integer numbers following the increasing ordering in the coordinates $x_1$ and
then in $x_2$. For the sake of clarity, we provide the following table that recapitulates the
domains and the corresponding expression of $F$. Since the map $F$ is continuous (but not $C^1$),
the expressions coincide on the domain boundaries. Accordingly, (when stability is not considered)
boundaries can be regarded as simultaneously belonging the adjacent domains.

\begin{figure}
\centerline{
\hbox{\psfig{figure=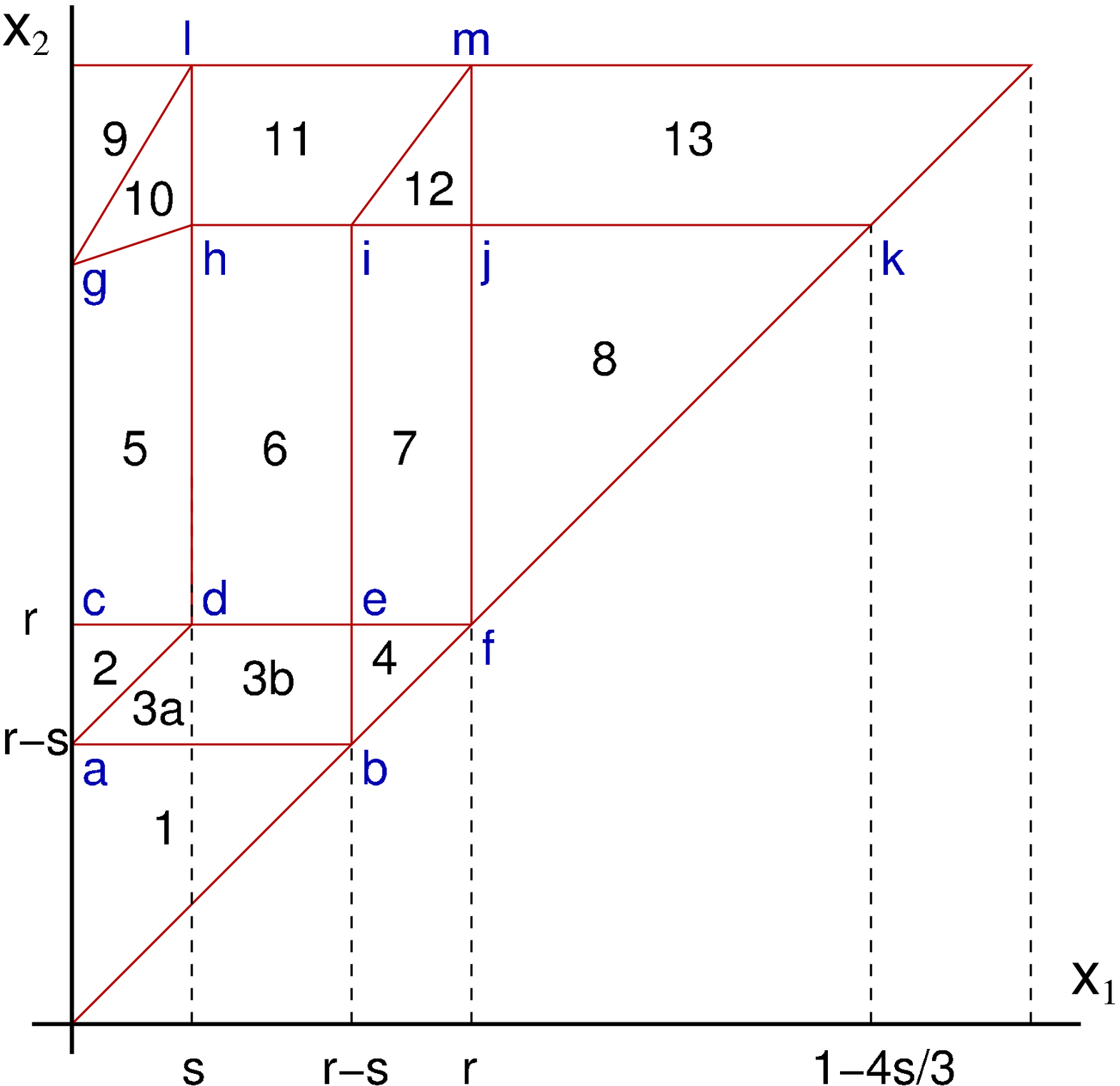,width=6cm,height=6cm} \hspace*{.4cm}
\psfig{figure=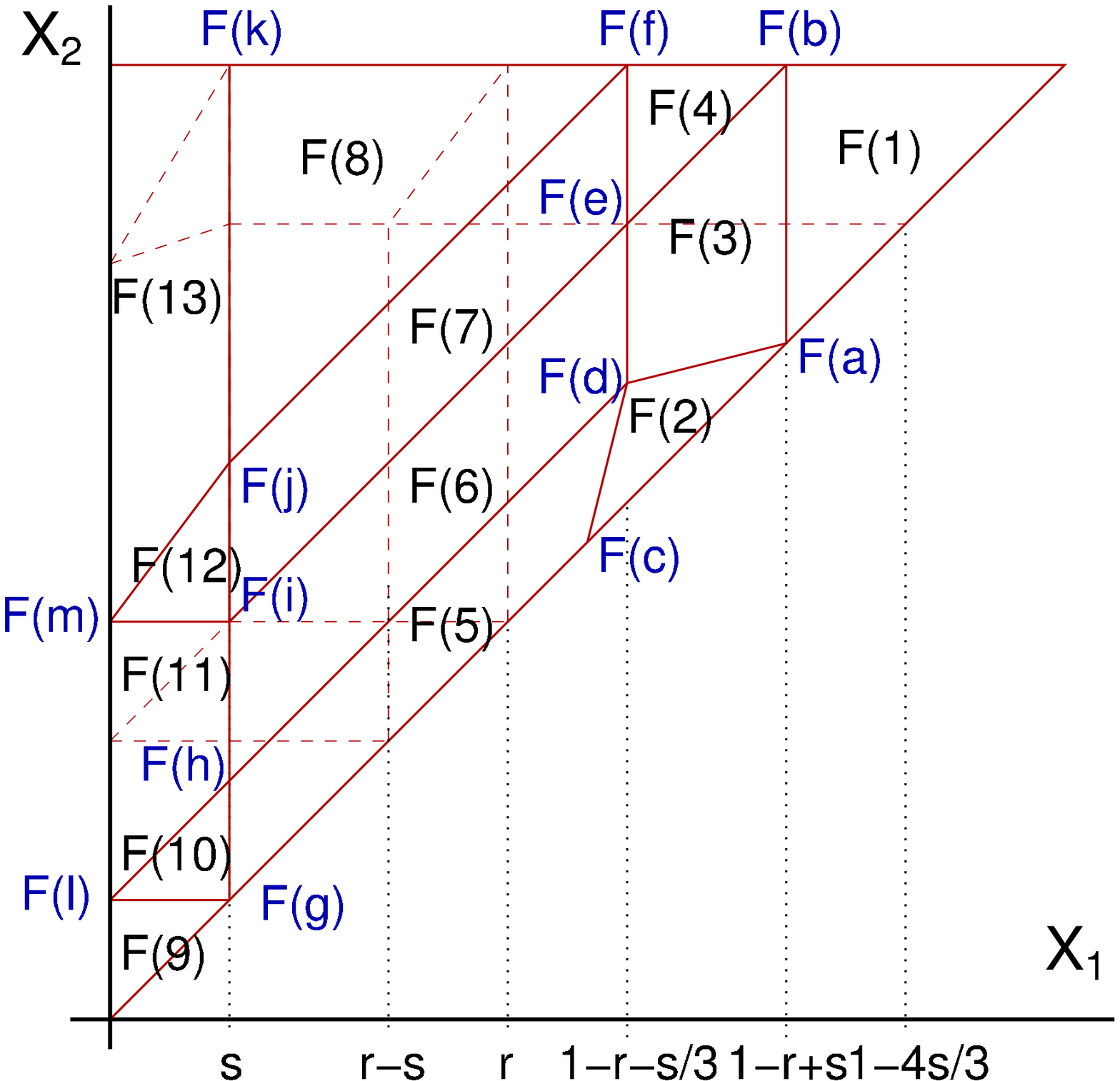,width=6cm,height=6cm}}
}
\caption{$k=3$. Partition of the triangle $0\leq x_1\leq x_2\leq 1$ and its image under $F$ for the pair $(r,s)=(\frac{5}{12},\frac{1}{8})$ chosen such that $2s< r\leq \frac{1}{2}-\frac{s}{3}$.}
\label{CASE1}
\end{figure}

{\small
$$
\begin{array}{|c|c|c|c|c|}
\hline
\text{Label}&\text{Condition\ on}\ x_1&\text{Condition\ on}\ x_2&x_0(t_1)&x_1(t_1)\\
\hline
1&0\leq x_1\leq x_2&0\leq x_2< r-s&1-x_2&1-x_2+x_1\\
\hline
3a&0\leq x_1\leq s&r-s\leq x_2< r-s+x_1&1-\frac{4x_2}{3}+\frac{r-s}{3}&1-\frac{4x_2}{3}+x_1+\frac{r-s}{3}\\
\hline
2&0\leq x_1\leq s&r-s+x_1\leq x_2<  r&1-\frac{5x_2-x_1}{3}+\frac{2(r-s)}{3}&1-\frac{5x_2-4x_1}{3}+\frac{2(r-s)}{3}\\
\hline
5&0\leq x_1\leq s&r\leq x_2<  1-\frac{5s}{3}+\frac{x_1}{3}&1-x_2+\frac{x_1}{3}-\frac{2s}{3}&1-x_2+\frac{4x_1}{3}-\frac{2s}{3}\\
\hline
10&0\leq x_1\leq s&1-\frac{5s-x_1}{3}\leq x_2<  1-\frac{5(s-x_1)}{3}&\frac{3(1-x_2)}{4}+\frac{x_1-s}{4}&\frac{3(1-x_2)}{4}+\frac{5x_1-s}{4}\\
\hline
9&0\leq x_1\leq s&1-\frac{5(s-x_1)}{3}\leq x_2\leq 1&\frac{3(1-x_2)}{5}&x_1+\frac{3(1-x_2)}{5}\\
\hline
3b&s<x_1\leq r-s&r-s\leq x_2<  r&1-\frac{4x_2}{3}+\frac{r-s}{3}&1-\frac{4x_2}{3}+x_1+\frac{r-s}{3}\\
\hline
6&s<x_1\leq r-s&r\leq x_2< 1-\frac{4s}{3}&1-x_2-\frac{s}{3}&1-x_2+x_1-\frac{s}{3}\\
\hline
11&s<x_1\leq r-\frac{3(1-x_2)}{4}&1-\frac{4s}{3}\leq x_2\leq 1&\frac{3(1-x_2)}{4}&x_1+\frac{3(1-x_2)}{4}\\
\hline
4&r-s<x_1\leq x_2&r-s\leq x_2<  r&1-\frac{4x_2}{3}+\frac{r-s}{3}&1-\frac{4(x_2-x_1)}{3}\\
\hline
7&r-s<x_1\leq r&r\leq x_2\leq 1-\frac{4s}{3}&1-x_2-\frac{s}{3}&1-x_2+\frac{4x_1}{3}-\frac{r}{3}\\
\hline
12&r-\frac{3(1-x_2)}{4}<x_1\leq r&1-\frac{4s}{3}\leq x_2\leq 1&\frac{3(1-x_2)}{4}&1-x_2+\frac{4x_1}{3}-\frac{r}{3}\\
\hline
8&r<x_1\leq x_2&r\leq x_2\leq 1-\frac{4s}{3}&1-x_2-\frac{s}{3}&1-x_2+x_1\\
\hline
13&r< x_1\leq 1&1-\frac{4}{3}s<x_2\leq 1&\frac{3(1-x_2)}{4}&1-x_2+x_1\\
\hline
\end{array}
$$
}

\subsection{The symbolic dynamics}
In order to obtain insights into the dynamics, we study the image of the triangle partition by using explicit
expressions in the table above. To that goal, we begin with a series of elementary observations.
\begin{itemize}
\item[$\bullet$] We already know from the analysis for arbitrary $k$ that the 3 corner points
are mapped one into another in a stable period-3 orbit. We also know that there must be another
3-periodic orbit on the edges.
\item[$\bullet$] In addition, there must be a fixed point inside the triangle.
\item[$\bullet$] $F(a)=F(0,r-s)=(1-r+s,1-r+s)$ and $F(b)=F(r-s,r-s)=(1-r+s,1)$. So the horizontal
segment $(a,b)$ is mapped into a vertical one with abscissa $1-r+s$. Depending on the location
of $r$ in the interval $\left[2s,1-\frac{5s}{3}\right]$, the quantity $1-r+s$ varies in the
interval $[s,1]$. Precisely, we have (one or several cases might not apply depending on $s$)
\[
\begin{array}{|c|c|}
\hline
\text{Conditions on }(r,s)&\text{Location of the segment} [F(a),F(b)]\\
\hline
2s\leq r\leq \frac{1+s}{2}&\text{beyond }r\\
\hline
\frac{1+s}{2}<r\leq \frac{1}{2}+s&\text{between }r-s\text{ and }r\\
\hline
\frac{1}{2}+s<r\leq 1-\frac{5s}{3}&\text{between }s\text{ and }r-s\\
\hline
\end{array}
\]
\item[$\bullet$] Similarly, the segment $(d,e,f)$ is mapped into a vertical one since we have $F(d)=F(s,r)=(1-r-\frac{s}{3},1-r+\frac{2s}{3})$, $F(e)=F(r-s,r)=(1-r-\frac{s}{3},1-\frac{4s}{3})$
and $F(f)=F(r,r)=(1-r-\frac{s}{3},1)$. The parameter dependent location of $1-r-\frac{s}{3}$ is
listed in the following table
\[
\begin{array}{|c|c|}
\hline
\text{Conditions on }(r,s)&\text{Location of the segment} [F(d),F(f)]\\
\hline
2s\leq r\leq \frac{1}{2}-\frac{s}{6}&\text{beyond }r\\
\hline
\frac{1}{2}-\frac{s}{6}<r\leq \frac{1}{2}+\frac{s}{3}&\text{between }r-s\text{ and }r\\
\hline
\frac{1}{2}+\frac{s}{3}<r\leq 1-\frac{5s}{3}&\text{between }s\text{ and }r-s\\
\hline
\end{array}
\]
In addition, we have $F(c)=F(0,r)=(1-r-\frac{2s}{3},1-r-\frac{2s}{3})$.
\item[$\bullet$] The segment $(g,h,i,j,k)$ is mapped into one with abscissa $s$ and in particular
 $F(g)=F(0,1-\frac{5s}{3})=(s,s)$, $F(h)=(s,2s)$, $F(i)=(s,r)=d$ $F(j)=(s,r+\frac{4s}{3})$ and $F(k)=(s,1)=l$.
\item[$\bullet$] We already know that the segment $(l,m)$ is mapped onto the vertical axis $x_1=0$.
Moreover $F(l)=(0,s)$ and $F(m)=(0,r)=c$.
\end{itemize}
The two previous properties imply that the regions in the upper vertical strip are mapped into
the following ones (see Figure \ref{CASE1})
\[
\begin{array}{l}
9,10\mapsto 1\\
11\mapsto 1\cup 2\cup 3a\\
12\mapsto 5\\
13\mapsto 5\cup 10\cup 9
\end{array}
\]
For the remaining regions, the situation depends on parameters and decomposes in various cases
according to the two previous tables. We consider here the case where all points from $a$ to
$f$ fall beyond the vertical line $x_1=r$, {\it i.e.} $2s< r\leq \frac{1}{2}-\frac{s}{3}$.
From Figure \ref{CASE1}, we get the following properties
\[
\begin{array}{l}
1,2,3,4\mapsto 8\cup 13\\
5,6\mapsto 1\cup 3\cup 4\cup 6\cup 7\cup 8\\
7\mapsto 6\cup 7\cup 8\cup 13\\
8\mapsto 6\cup 7\cup 11\cup 12\cup 13
\end{array}
\]

\subsection{Parameter dependence of fixed points and bifurcations}
The above properties suggest to consider possible fixed points in 6 and in 7. It turns out
more convenient to study the one in 7 and to follow its bifurcations in parameter space.

\begin{figure}
\begin{center}
\psfig{figure=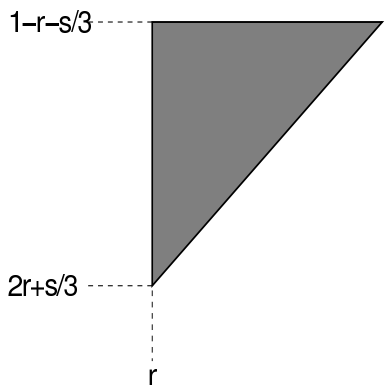,width=5cm,height=5cm}
\end{center}
\caption{Invariant triangle composed of the neutral fixed point and neutral period 3 orbits lying in 8.
The 3 corners are points of the 3-periodic orbit with code $7\mapsto 7\mapsto 8\mapsto 7$. Notice that
the left vertical edge (belongs to 8 and) is part of the boundary between 7 and 8.}
\label{NEUTRALTRIANGLE}
\end{figure}

\underline{Fixed point in 7:} The (unique) fixed point in 7, namely
$(\frac{3 + r - 2 s}{10},\frac{21 - 3 r - 4 s}{30})$ turns out to be a source
(with associated complex eigenvalues for $F^3$ - the return map under study) which exists provided that
\[
\frac{3-2s}{9}\leq r\leq \frac{3+8s}{9}
\]
At the lower boundary $r=\frac{3-2s}{9}$ of this domain, the first coordinate meets the vertical
boundary with 8. This suggests that $(\frac{3 + r - 2 s}{10},\frac{21 - 3 r - 4 s}{30})$ might be
continued as a fixed point in 8 for $r<\frac{3-2s}{9}$. This is indeed the case.
\begin{itemize}
\item[$\bullet$] \underline{Fixed point in 8:} Easy calculations conclude that the coordinates are
given by $(\frac{3-2s}{9},\frac{6-s}{9})$. Moreover, the fixed point is neutral (with a double
eigenvalue of $F^3$ equal to 1) and it exists provided that
\[
2s\leq r\leq \frac{3-2s}{9}
\]
For $r=\frac{3-2s}{9}$ this point coincides with the fixed point in 7 and we expect a kind of
pitchfork bifurcation which would create a 3-periodic orbit with code
$7\mapsto 7\mapsto 8\mapsto 7$.
\item[$\bullet$] \underline{Period-3 orbit $7\mapsto 7\mapsto 8\mapsto 7$:} Calculations show that
this (unique) orbit has coordinates:
\[
(r,1-r-\frac{s}{3})\mapsto (r,2r+\frac{s}{3})\mapsto (1-2r-\frac{2s}{3},1-r-\frac{s}{3})\mapsto (r,1-r-\frac{s}{3})
\]
and indeed emerges from the fixed point in 7 at $r= \frac{3-2s}{9}$. It persists for all
lower values of $r$ down to $2s$ and is expanding with double real eigenvalues.

\item[$\bullet$] \underline{Period-3 orbits $8\mapsto 8\mapsto 8\mapsto 8$:} Topological considerations
suggest that other orbits should come into play at $r= \frac{3-2s}{9}$. Indeed, calculations confirm
that a two-parameter family of 3-periodic orbits lying in 8 also emerges from the fixed point in 7
when $r$ crosses $\frac{3-2s}{9}$ downward. All orbits are neutral with double eigenvalue 1 and
their coordinates have the following expressions
\[
(x_1,x_2)\mapsto (1-x_2-\frac{s}{3},1-x_2+x_1)\mapsto (x_2-x_1-\frac{s}{3},1-\frac{s}{3}-x_1)\mapsto (x_1,x_2)
\]
where $x_1$ is arbitrary in the interval $\left(r,x_2-r-\frac{s}{3}\right)$, $x_2$ is arbitrary
in $\left(2r+\frac{s}{3},1-r-\frac{s}{3}\right)$ and $r$ is arbitrary in $\left[2s,\frac{3-2s}{9}\right]$.
The orbits can be viewed as rotating around the fixed point $(\frac{3-2s}{9},\frac{6-s}{9})$ which is
included in the family ({\it i.e.}\ for $(x_1,x_2)=(\frac{3-2s}{9},\frac{6-s}{9})$, the orbit actually
reduces to a fixed point). The closure of the set of orbit coordinates forms a triangle whose corners
are the coordinates of the orbit $7\mapsto 7\mapsto 8\mapsto 7$ (see Figure \ref{NEUTRALTRIANGLE}).
\end{itemize}

\noindent
At the upper boundary $r=\frac{3+8s}{9}$ of the ``7''-fixed point's existence domain, its first coordinate
meets the vertical boundary with 6. We then consider the fixed point in 6.
\begin{itemize}
\item[$\bullet$] \underline{Fixed point in 6:} Its coordinates are given by
$(\frac{3-s}{9},\frac{2(3-s)}{9})$ and this neutral fixed point exists iff
\[
\frac{3+8s}{9}\leq r\leq \frac{2(3-s)}{9}
\]
The bifurcation scenario at $r=\frac{3+8s}{9}$ is similar to the one above at the lower boundary
$r= \frac{3-2s}{9}$.
\item[$\bullet$] \underline{Period-3 orbit $7\mapsto 7\mapsto 6\mapsto 7$:} It turns out that
this unique orbit emerges from the fixed point 7 at $r=\frac{3+8s}{9}$ and exists up to $r=\frac{3+2s}{6}$
(or up to $1-\frac{5s}{3}$ whichever is smaller). It is expanding with double real eigenvalues.
At $r=\frac{3+2s}{6}$ both the component in 6 and its successor cross the horizontal line $x_2=r$
and the second component in 7 cross the upper domain boundary $x_2=1-\frac{4s}{3}$.

However, the value $r=\frac{3+2s}{6}$ does not involve any existence condition related to the fixed
point in 6. As shown below, the reason is that the periodic orbit can be continued  up to
$\frac{2(3-s)}{9}$ provided that we consider the appropriate sequence of symbols.

\item[$\bullet$] \underline{Period-3 orbit $6\mapsto 3\mapsto 3\mapsto 6$:} This orbit continues the
orbit $7\mapsto 7\mapsto 6\mapsto 7$ when $r>\frac{3+2s}{6}$. Indeed, it exists provided that
\[
\frac{3+2s}{6}\leq r\leq \frac{2(3-s)}{9}
\]
It is expanding with 2 real eigenvalues (the same as those associated with $7\mapsto 7\mapsto 6\mapsto 7$)
and coincide with $7\mapsto 7\mapsto 6\mapsto 7$ at $r=\frac{3+2s}{6}$.

\item[$\bullet$] \underline{Period-3 orbits $6\mapsto 6\mapsto 6\mapsto 6$:} Similarly to as before,
a triangle of neutral 3-periodic orbits is created from the fixed point 7 at $r=\frac{3+8s}{9}$. The
coordinates have the following expression
\[
(x_1,x_2)\mapsto (1-x_2-\frac{s}{3},1-x_2+x_1-\frac{s}{3})\mapsto (x_2-x_1,1-x_1-\frac{s}{3})\mapsto (x_1,x_2)
\]
and these orbits exist in the same parameter domain as the fixed point in 6. When $r\in \left[\frac{3+8s}{9},\frac{3+2s}{6}\right]$ the restrictions on coordinates are
\[
1-2r+\frac{5s}{3}< x_1< r-s\quad\text{and}\quad 1-r+\frac{2s}{3}< x_2< x_1+r-s
\]
and they are
\[
2r+\frac{s}{3}-1< x_1< 1-r-\frac{s}{3}\quad\text{and}\quad r< x_2< x_1+1-r-\frac{s}{3}
\]
when $r\in \left[\frac{3+2s}{6},\frac{2(3-s)}{9}\right]$.
\end{itemize}

At the boundary $r=\frac{2(3-s)}{9}$, the expanding periodic orbit and the family of neutral
orbits meet with the neutral fixed point ``6'' in an inverse pitchfork bifurcation to create the
fixed point in 3.
\begin{itemize}
\item[$\bullet$] \underline{Fixed point in 3:} The fixed point $(\frac{3+r-s}{11},\frac{2(3-r-s)}{11})\in 3$
is expanding with double real eigenvalues and exists in the domain
\[
\frac{2(3-s)}{9}\leq r\leq \frac{2}{3}+s
\]
\end{itemize}

At the upper boundary $r=\frac{2}{3}+s$ of its existence domain, the fixed point "3" meets the region 1.
\begin{itemize}
\item[$\bullet$] \underline{Fixed point in 1:} The fixed point $(\frac{1}{3},\frac{2}{3})$ is
neutral with double real eigenvalue 1 and exists in the domain
\[
\frac{2}{3}+s\leq r< 1
\]

\item[$\bullet$] \underline{Period-3 orbit $3\mapsto 3\mapsto 1\mapsto 3$:} Same parameter domain.
Expression of coordinates
\[
(2(r-s)-1,r-s)\mapsto (1-r+s,r-s)\mapsto (1-r+s,2(1-r+s))\mapsto (2(r-s)-1,r-s)
\]
Expanding with double real eigenvalue. Emerges from the fixed point 3.

\item[$\bullet$] \underline{Period-3 orbit $1\mapsto 1\mapsto 1\mapsto 1$:} Same parameter domain.
Triangle of neutral 3-periodic orbits
\[
(x_1,x_2)\mapsto (1-x_2,1-x_2+x_1)\mapsto (x_2-x_1,1-x_1)\mapsto (x_1,x_2)
\]
where $x_1$ is arbitrary in $\left(1-r+s,2(r-s)-1\right)$ and $x_2\in \left(1-r+s+x_1,r-s\right)$. As
before, the corners coincides with the period-3 orbit $3\mapsto 3\mapsto 1\mapsto 3$.
\end{itemize}

\begin{figure}
\begin{center}
\epsfig{figure=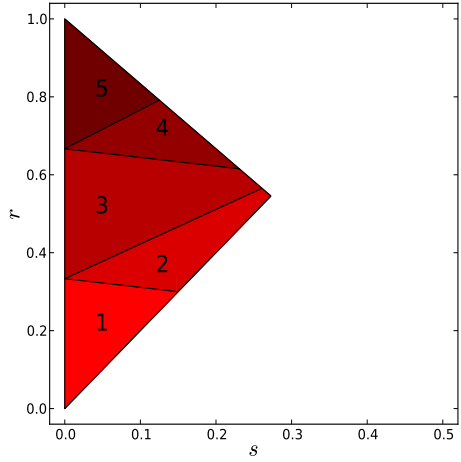,width=7cm,height=7cm}
\end{center}
\caption{$k=3$. Parameters domains of existence of fixed points and
related periodic orbits analyzed in this section. In regions (2) and (4) the unique fixed point
is unstable. In regions (1), (3) and (5) the fixed point is neutral and sits inside a triangle of
neutral period 3 points.}
\label{DOMAINS}
\end{figure}

To recapitulate, the return map $F^3$ possesses the following orbits depending on parameters in
the considered regions (see Figure \ref{DOMAINS})
\begin{itemize}
\item[(1)] $2s\leq r\leq \frac{3-2s}{9}$. three sources (two of them belong to 7, the other one
lies in 8) and a two-parameter family of neutral fixed points lying inside the triangle whose corners
are the 3 sources.
\item[(2)] $\frac{3-2s}{9}<r\leq \frac{3+8s}{9}$. One source in 7.
\item[(3)] $\frac{3+8s}{9}<r\leq \frac{2(3-s)}{9}$. Similarly to as in (1); namely three sources (two in 7,
one 6 if $r \leq \frac{3+2s}{6}$, and two sources belong to 3 and one is in 6 for
$r \geq \frac{3+2s}{6}$)
and a triangle of neutral fixed points in 6.
\item[(4)] $\frac{2(3-s)}{9}<r\leq \frac{2}{3}+s$. Similarly to as in (2); namely a single source in 1.
\item[(5)] $\frac{2}{3}+s<r< 1-\frac{5s}{3}$. Similarly to as in (1); namely three sources
(two lie in 3 and 1 lies in 1) and a triangle of neutral fixed points in 1. 
\end{itemize}

Note that in case (5) we have $r-s > \frac{2}{3}$ and this corresponds to the case that the three
clusters in the cyclic solution (which has initial conditions $(0,1/3,2/3)$ do not interact, {\sl i.e.}\ 
$x_0$ leaves $S$ before $x_2$ enters $R$ and no feedback is experienced.

In the cases (1) and (3) the clusters in the cyclic solution experience feedback, but in a non-essential way,
by which we mean that small perturbations do not lead either to further separation or contraction between
the clusters. For instance in case (3) $x_2$ begins in $R$ and $x_1$ begins in between $S$ and $R$ and
$x_0$ leaves $S$ before either of these states changes. In case (1) both $x_1$ and $x_2$ begin in
$R$, but $x_0$ leaves $S$ before either leaves $R$.

A python script that will produce movies for the dynamics for arbitrary $0<r\le s <1$ can be found
at: \verb&http://oak.cats.ohiou.edu/~rb301008/research.html&.

\subsection{Additional orbits}
Besides orbits bifurcating with fixed points, based on numerics and on properties of the symbolic
dynamics, additional orbits of $F$ exist depending on parameters. Their existence domains do not
coincide with those listed above.

\begin{itemize}
\item[$\bullet$] The most important orbits are the one-parameter family of period-3 neutral-stable
orbits lying on the edges (and with code $1\mapsto 1\mapsto 11\mapsto 1$)
\[
(0,x_2)\mapsto (1-x_2,1-x_2)\mapsto (x_2,1)
\]
These orbits exist for arbitrary $x_2\in \left(1-r+s,r-s\right)$ provided that
$\frac{1}{2}+s\leq r<1-\frac{5s}{3}$. They are stable with respect to transverse perturbations.
The coordinates form 3 open intervals, each included in one the edges. Moreover, the interval
boundaries respectively form a period-3 hyperbolic orbit (code $2\mapsto 1\mapsto 11\mapsto 2$)
and a period-3 neutral-unstable orbit (code $1\mapsto 4\mapsto 11\mapsto 1$) which exist in the
same parameter domains and which merge for $r=\frac{1}{2}+s$. Orbits in the family correspond to
2 clusters solutions in the original system with one cluster composed of two clusters (and the
two clusters being certainly isolated one from each other).
\item[$\bullet$] Interestingly, the bifurcation at $r=\frac{1}{2}+s$ generates two distinct periodic
orbits with identical code ($2\mapsto 4\mapsto 11\mapsto 2$). Both orbits are neutral-unstable and
the first has coordinates
\[
(\frac{3-6(r-s)}{4},\frac{3-2(r-s)}{4})\mapsto (r-s,\frac{3-2(r-s)}{4})\mapsto (r-s,2(r-s))
\]
with the first point being in 2. It exists under the condition $\frac{1}{2}+\frac{s}{3}\leq r\leq \frac{1}{2}+s$.

The other solution is actually a one-parameter family of periodic orbits
(which forms a segment and) whose component in 11 is $\left(\frac{5(s-r)+3x_2}{11},x_2\right)$
where $x_2$ can be chosen arbitrary with the condition
\[
\frac{3}{4}(1-x_2)\leq \min\left\{\frac{3-6(r-s)}{5},\frac{6r+5s-3}{7}\right\}
\]
The family exists under the condition $\frac{1}{2}-\frac{5s}{6}\leq r\leq \frac{1}{2}+s$.
\end{itemize}
The bifurcations taking place at $r=\frac{1}{2}+\frac{s}{3}$ and $r=\frac{1}{2}-\frac{5s}{6}$
are unclear. The fact that the one-parameter family reaches the boundary with 5 in the former
case suggests to investigate orbits with code $5\mapsto 4\mapsto 11\mapsto 5$. Surprisingly,
the analysis concludes that such an orbit exits only if $r=\frac{1}{2}+\frac{s}{3}$, a value
that is unrelated to the one-parameter family existence condition.

Two additional orbits have been found on the edge when $2s<r\leq \frac{1}{2}+\frac{s}{12}$.
One is hyperbolic with code $2\mapsto 8\mapsto 11\mapsto 2$ and coordinates
\[
(0,r-\frac{s}{2})\mapsto (1-r+\frac{s}{6},1-r+\frac{s}{6})\mapsto (r-\frac{s}{2},1)\mapsto (0,r-\frac{s}{2}).
\]
The other one is neutral-unstable with code $2\mapsto 4\mapsto 11\mapsto 2$ and coordinates
\[
(0,\frac{3+5(r-s)}{11})\mapsto (1-\frac{5+r-s}{11},1-\frac{5+r-s}{11})\mapsto (\frac{3+5(r-s)}{11},1)\mapsto (0,\frac{3+5(r-s)}{11}).
\]
The two orbits merge for $r=\frac{1}{2}+\frac{s}{12}$ in a seemingly saddle-node bifurcation.


\section{Discussion}

In this note we have investigated only a portion of the possible parameter space $0<s \le r <1$.
It should be clear to the reader by this point that further investigations into other subsets of the 
parameters in this fashion are possible, but perhaps prohibitively time-consuming.
We see that such studies are likely to also prove unprofitable, since numerical simulations
show that no other types of dynamics occur other than those described here.
One can download a python script that to investigate the dynamics for arbitrary $0<r\le s <1$ 
at: \verb&http://oak.cats.ohiou.edu/~rb301008/research.html&.

The dynamics we have observed for these parameter sets closely resembles the dynamics of $k=2$ cluster
systems analyzed in \cite{YFBMB}. The $k=3$ fixed point of $F$ with positive feedback, like that for
$k=2$ is either:
\begin{itemize}
\item isolated and locally unstable, or,
\item contained in a neutrally stable set of period $k$ points.
\end{itemize}
In the later case the edges of the neutrally stable set are unstable. This case exists if either the
three clusters are isolated from each other, or, if they interact in a non-essential way.
In both cases the orbits of all other interior points are asymptotic to the boundary. Thus $k$ cyclic 
solutions for either $k=2$ or $k=3$ are practically
unstable in the sense the arbitrarily small perturbations may lead to loss of stability and eventual
 merger of clusters. Since the single cluster solution (synchronization) is the only solution that
 is asymptotically stable, it would seem to be the most likely to be observed in application if
 the feedback is similar to the form we propose and is positive.

{\bf Acknowledgments:}\\
B.F.~thanks the Courant Institute (NYU) for hospitality. He~was supported by CNRS
and by the EU Marie Curie fellowship PIOF-GA-2009-235741.
T.Y.\ and this work were supported by the NIH-NIGMS grant R01GM090207.

\end{document}